\documentclass{amsart}
\usepackage{amsmath}
\usepackage{amssymb}
\usepackage{amsthm}
\usepackage{color}

\newcommand{\bC}{\mathbb{C}}
\newcommand{\bH}{\mathbb{H}}
\newcommand{\bR}{\mathbb{R}}
\newcommand{\bZ}{\mathbb{Z}}

\newcommand{\mF}{\mathcal{F}}
\newcommand{\lv}{\lvert}
\newcommand{\rv}{\rvert}
\newcommand{\lp}{\langle}
\newcommand{\rp}{\rangle}

\newcommand{\delbar}{\overline{\partial}}
\newcommand{\hA}{\widehat{A}}
\newcommand{\htheta}{\widehat{\theta}}
\newcommand{\hnabla}{\widehat{\nabla}}
\DeclareMathOperator{\Real}{Re}

\newcommand{\suchthat}{\mathrel{}\middle|\mathrel{}}

\newtheorem{thm}{Theorem}

\newtheorem{pquestion}{Question}

\theoremstyle{remark}

\numberwithin{equation}{section}

\begin{document}

\title[CR Lichnerowicz conjecture with mixed signature]{On the Lichnerowicz conjecture for CR manifolds with mixed signature}
\author{Jeffrey S.\ Case}
\address{109 McAllister Building \\ Penn State University \\ University Park, PA 16802}
\email{jscase@psu.edu}
\author{Sean N.\ Curry}
\address{Department of Mathematics\\ University of California, San Diego\\ 9500 Gilman Drive \\ La Jolla, CA 92093-0112, USA}
\email{sncurry@math.ucsd.edu}
\author{Vladimir S.\ Matveev}
\address{Institut f\"ur Mathematik\\ Fakult\"at f\"ur Mathematik und Informatik \\ Friedrich-Schiller-Universit\"at Jena \\ 07737 Jena, Germany}
\email{vladimir.matveev@uni-jena.de}
\begin{abstract}
 We construct examples of nondegenerate CR manifolds with Levi form of signature $(p,q)$, $2\leq p\leq q$, which are compact, not locally CR flat, and admit essential CR vector fields.  We also construct an example of a noncompact nondegenerate CR manifold with signature $(1,n-1)$ which is not locally CR flat and admits an essential CR vector fields.  These provide counterexamples to the analogue of the Lichnerowicz conjecture for CR manifolds with mixed signature.
\end{abstract}
\maketitle

\section{Introduction}
\label{sec:intro}

A \emph{conformal vector field} on a Riemannian manifold $(M^n,g)$ is a vector field $X$ whose flow acts by conformal diffeomorphisms.  Such a vector field is \emph{essential} if there is no smooth metric $\hat g=e^{2\Upsilon}g$ conformal to $g$ which is preserved by the flow of $X$.  For example, the infinitesimal generator of Euclidean dilations is an essential conformal vector field on Euclidean space, and the infinitesimal generator of the conjugation of Euclidean dilations with stereographic projection is an essential conformal vector field on the round sphere.  A famous conjecture attributed to Lichnerowicz hypothesized that these are the only Riemannian manifolds, up to conformal diffeomorphism, which admit essential conformal vector fields.  This was verified by Obata~\cite{Obata1971} and Ferrand~\cite{Ferrand1971,Ferrand1996}.

A \emph{CR vector field} on a pseudohermitian manifold $(M^{2n+1},J,\theta)$ is a vector field $X$ which preserves the contact distribution $H:=\ker\theta$ and the complex structure $J$ on $H$.  Such a vector field is \emph{essential} if there is no contact form $\hat\theta=e^\Upsilon\theta$ which is preserved by the flow of $X$.  For example, the infinitesimal generator of Heisenberg dilations on $\bH^n$ and the infinitesimal generator of the conjugation of Heisenberg dilations with the Cayley transform are essential CR vector fields on $\bH^n$ and the round CR $(2n+1)$-sphere, respectively.  Schoen~\cite{Schoen1996} gave an alternative proof of the Lichnerowicz Conjecture which also implied that the only strictly pseudoconvex manifolds, up to CR diffeomorphism, which admit essential CR vector fields are the Heisenberg groups and the round CR spheres.

It is an interesting question to ask whether there are similar classifications of conformal and CR manifolds of mixed signature which admit essential vector fields.  This question was posed explicitly by D'Ambra and Gromov~\cite{DAmbraGromov1991} in the conformal category for closed manifolds, the reason being that noncompact counterexamples of Lorentzian signature were known.  The CR analogue of their question is as follows:

\begin{pquestion}
 \label{q:lichnerowicz}
 Let $(M^{2n+1},J,\theta)$ be a closed nondegenerate pseudohermitian manifold which admits an essential CR vector field.  Is $(M^{2n+1},J,\theta)$ locally CR flat?
\end{pquestion}

In this note we give a negative answer to Question~\ref{q:lichnerowicz} in the case when the Levi form has signature $(p,q)$, $p,q\geq2$.  Indeed, we construct infinitely many CR inequivalent examples.  We also show that there are noncompact nondegenerate pseudohermitian manifolds with Levi form of signature $(1,n-1)$ which admit essential CR vector fields.  Our counterexamples to Question~\ref{q:lichnerowicz} are heavily inspired by corresponding counterexamples in the conformal case given by Frances~\cite{Frances2015}.

A negative answer to Question~\ref{q:lichnerowicz} is readily constructed using Chern--Moser normal form~\cite{ChernMoser1974}, as we describe explicitly for the case of signature $(2,2)$.  To that end, let $w\in\bC$ and $z\in\bC^4$; we write $z=(z_1,\dotsc,z_4)$.  Set
\[ h(z,\bar z) = z_1\bar z_2 + z_2\bar z_1 + z_3\bar z_4 + z_4\bar z_3, \]
so that $h$ is a nondegenerate hermitian form of signature $(2,2)$.  Define the real function $r\colon\bC^5\to\bR$ by
\begin{equation}
 \label{eqn:r}
 r(w,z) := \frac{1}{2i}(w-\bar w) - h(z,\bar z) - \lv z_1\rv^4 .
\end{equation}
Since $z\mapsto \lv z_1\rv^4$ is trace-free with respect to $h$, the work of Chern and Moser~\cite{ChernMoser1974} implies that the level set $M:=r^{-1}(0)$ is a nondegenerate CR manifold with Levi form of signature $(2,2)$ which is not locally CR flat.  Next, note that for any $\alpha,\beta\in\bR$, the linear action $\Gamma_{\alpha,\beta}\colon\bC^5\to\bC^5$ defined by
\[ \Gamma_{\alpha,\beta}(w,z_1,z_2,z_3,z_4) := \left( e^{4\beta}w, e^\beta z_1, e^{3\beta} z_2, e^{4\beta-\alpha}z_3, e^\alpha z_4 \right) \]
fixes $M$ and acts by homotheties of the contact form $\theta := \mathrm{Re}(i\partial r)$.  In particular, with $4\beta<\alpha<0$ fixed, $\Gamma_{\alpha,\beta}$ generates a free $\bZ$-action with respect to which $M\setminus\{0\}$ has compact quotient $M_{\alpha,\beta}$; note that $M_{\alpha,\beta}$ is diffeomorphic to $S^1\times S^8$.  Since the one-parameter family $\phi_t:=\Gamma_{0,-t}$ commutes with $\Gamma_{\alpha,\beta}$, and hence descends to a CR flow on $M_{\alpha,\beta}$, the infinitesimal generator $X$ of $\phi_t$ is necessarily a CR vector field.  Moreover, the (real) two-dimensional torus obtained as the image of the fixed set $\{w=z_1=z_2=z_3=0\}\subset M\setminus\{0\}$ is an attractor of the flow of $\phi_t$ on $M_{\alpha,\beta}$.  In particular, $\phi_t$ cannot even preserve a Borel measure which is nonzero on open sets, let alone a contact form (cf.\ \cite{Frances2015}).

The above construction provides infinitely many geometrically-distinct counterexamples to Question~\ref{q:lichnerowicz}.  To see this, define $\gamma\colon\bC\to M$ by
\[ \gamma(z) = (0,0,z,0,0) . \]
The map $\gamma$ is compatible with the complex structure on $M$, meaning $\gamma_\ast\partial_z\in T^{1,0}M$ and $\gamma_\ast\partial_{\bar z}\in T^{0,1}M$.  The image of $\gamma$ is the unique integral surface of the image of the Chern tensor, regarded as a $(1,3)$-tensor, with the property that it is $\Gamma_{\alpha,\beta}$-invariant.  It turns out that $\gamma$ can be given a CR invariant projective structure which is analogous to the conformally invariant projective structure on null geodesics~\cite{Markowitz1981}.  Using this structure, one can show that the quotients $M_{\alpha,\beta}$ and $M_{\tilde\alpha,\tilde\beta}$ are CR inequivalent when $\beta\not=\tilde\beta$.

To negatively answer Question~\ref{q:lichnerowicz} in the case when the Levi form has signature $(p,q)$, $p,q\geq 2$, we replace $h$ by the nondegenerate hermitian form
\[ h(z,\bar z) = z_1\bar z_2 + z_2\bar z_1 + z_3\bar z_4 + z_4\bar z_3 + \sum_{j=5}^{p+2} z_j\bar z_j - \sum_{k=p+3}^{p+q} z_k\bar z_k \]
of signature $(p,q)$ on $\bC^{p+q}$ and modify the linear action $\Gamma_{\alpha,\beta}$ so that $z_j\mapsto e^{2\beta} z_j$ for $5\leq j\leq p+q$.  With these modifications, the rest of the above argument carries through.

The essential ingredient in the above argument is the existence of a signature $(p,q)$ hypersurface in Chern-Moser normal form which is not locally CR flat and admits two commuting non-isometric homotheties, one of which induces a cocompact free $\bZ$-action and one of which has a nontrivial fixed set.  The former homothety yields a compact quotient, while the infinitesimal generator of the latter homothety is an essential CR vector field.  While our construction does not yield a negative answer to Question~\ref{q:lichnerowicz} with Levi form of signature $(1,n-1)$, it does yield a negative answer for noncompact manifolds: By considering the defining function and group actions on $\bC^3$ obtained by omitting the variables $z_3$ and $z_4$, one obtains a non-isometric homothety which fixes the origin, and thus its infinitesimal generator is an essential CR vector field.

\begin{thm}
 \label{thm:main}
 For any $p,q\geq2$, there exist infinitely many pairwise inequivalent non-flat CR structures of signature $(p,q)$ on $S^1\times S^{2p+2q}$ which admit an essential CR vector field.  Moreover, there is a non-flat CR structure of signature $(1,n-1)$ on $\bR\times\bC^n$ which admits an essential CR vector field.
\end{thm}

We conclude by pointing out that our example $M:=r^{-1}(0)$ for the defining function~\eqref{eqn:r} has appeared elsewhere, namely as an example of a CR manifold with submaximal symmetry dimension~\cite{Kruglikov2016}.  This gives a further parallel to the conformal counterexamples found by Frances; indeed, results of Kruglikov and The~\cite[\S 5.1.1]{KruglikovThe2017} imply that Frances' metric~\cite[(1)]{Frances2015} has submaximal conformal symmetry dimension.

The rest of this note is organized as follows: In Section~\ref{sec:markowitz}, we describe an invariant projective structure on null geodesics in nondegenerate CR manifolds.  In Section~\ref{sec:details}, we use this structure to complete the proof of Theorem~\ref{thm:main}.

\subsection*{Acknowledgements} The work described in this note was done during and shortly after the workshop {\it Analysis and geometry on pseudohermitian manifolds} held at the American Institute of Mathematics (AIM) and organized by Sorin Dragomir, Howard Jacobowitz, and Paul Yang in November 2017. The style of AIM workshops is to break into small groups in the afternoon to work on a specific problem of interest to participants.  This note is the result of the discussions of one such group.  We thank AIM for a perfect working environment and atmosphere and for financial support, and other participants of the workshop, in particular Stephen McKeown, Howard Jacobowitz, Paul Yang and Sagun Chanillo for participation in our discussions and useful comments.

JSC was supported by a grant from the Simons Foundation (Grant No.\ 524601).

\section{An invariant projective structure on null geodesics}
\label{sec:markowitz}

In this section we show that if a CR manifold $(M^{2n+1},J)$, $n\geq2$, locally admits pseudo-Einstein contact forms, then any complex null geodesic in $M$ carries a CR invariant projective parameterization.  The former assumption means that for any point $p\in M$, there is a neighborhood of $p$ on which a contact form $\theta$ for $(M^{2n+1},J)$ can be chosen such that the pseudohermitian Ricci curvature $R_{\alpha\bar\beta}$ and the Levi form $L_\theta$ of $\theta$ are proportional.  This assumption allows us to easily adapt the construction by Markowitz~\cite{Markowitz1981} of a conformally invariant projective structure on null geodesics.

A \emph{complex null curve} in $M$ is a complex curve $\gamma\colon U\to M$, $U\subseteq\bC$ open, for which
\begin{enumerate}
 \item $\gamma^\prime:=\gamma_\ast\partial_z$ takes its values in $T^{1,0}$, and
 \item $\gamma^\prime$ is null with respect to $L_\theta$.
\end{enumerate}
Thus $\bar\gamma^\prime:=\gamma_\ast\delbar_z$ is also null.  Note that if $w\colon\bC\to\bC$ is a holomorphic change of variables, then $\gamma\circ w$ is also a complex null curve.

Suppose that $Z\in\Gamma\bigl(T^{1,0}\bigr)$ is a section of $T^{1,0}$ such that $[Z,\bar Z]=0$.  This implies that $Z$ is null with respect to $L_\theta$ and that the distribution $D:=\Real\lp Z\rp$ is an integrable rank two distribution.  A leaf of $D$ can be given a parameterization which makes it a complex null curve.

A \emph{complex null geodesic} of $(M,\theta)$ is a complex null curve $\gamma$ such that
\begin{equation}
 \label{eqn:geodesic}
 \nabla_{\gamma^\prime}\gamma^\prime = u\gamma^\prime, \qquad \nabla_{\bar\gamma^\prime}\gamma^\prime = 0
\end{equation}
for some holomorphic function $u$ along $\gamma$; i.e.\ we require that~\eqref{eqn:geodesic} holds for some complex-valued function $u$ defined along $\gamma$ which satisfies $\bar\gamma^\prime u = 0$.  This notion is clearly invariant under holomorphic change of variables.  An \emph{affine parameterization of $\gamma$ with respect to $\theta$} is a reparameterization $w=w(z)$ such that
\[ \nabla_{\dot\gamma}\dot\gamma = 0, \]
where dots denote differentiation with respect to $w$.  One readily checks that this occurs if and only if
\[ \frac{w^{\prime\prime}}{w^\prime} = u, \]
where primes denote differentiation with respect to $z$.  Since $u$ is holomorphic along $\gamma$, such a parameterization always exists locally.

A key point is that complex null geodesics are \emph{secondary invariants}; i.e.\ if $\gamma$ is a complex null geodesic with respect to a pseudo-Einstein contact form $\theta$ for $(M,J)$, then $\gamma$ is a complex null geodesic with respect to any pseudo-Einstein contact form for $(M,J)$.  To see this, recall that Lee~\cite{Lee1986} showed that for any vector fields $Z,W\in\Gamma\bigl(T^{1,0}\bigr)$, if $\htheta=e^\Upsilon\theta$, then
\begin{align*}
 \hat\nabla_W Z & = \nabla_W Z + W(\Upsilon)\,Z + Z(\Upsilon)\,W , \\
 \hat\nabla_{\bar W}Z & = \nabla_{\bar W}Z - h(Z, \bar W)\nabla^{1,0}\Upsilon,
\end{align*}
where $h$ is the Levi form associated to $\theta$.  In particular, if $\gamma\colon\bC\to M$ is a complex null curve, then
\begin{align*}
 \hnabla_{\gamma^\prime}\gamma^\prime & = \bigl( u + 2\Upsilon^\prime \bigr) \gamma^\prime, \\
 \hnabla_{\bar\gamma^\prime}\gamma^\prime & = 0 .
\end{align*}
If $\htheta$ and $\theta$ are pseudo-Einstein, then $\Upsilon$ is CR pluriharmonic~\cite{Lee1988}, and hence
\[ \bar\gamma^\prime \Upsilon^\prime = (\gamma^\prime)^\alpha(\gamma^\prime)^{\bar\beta}\Upsilon_{\alpha\bar\beta} = (\gamma^\prime)^\alpha(\gamma^\prime)^{\bar\beta}\Upsilon_\rho{}^\rho h_{\alpha\bar\beta}, \]
where we use abstract index notation to denote by $\Upsilon_\alpha$ the $(1,0)$-part of $d\Upsilon$ and by $\Upsilon_{\alpha\bar\beta}$ the $(0,1)$-part of the covariant derivative of $\Upsilon_\alpha$ with respect to the Tanaka--Webster connection.  Since $\gamma$ is null, we see that $\Upsilon^\prime$ is holomorphic.  This verifies that complex null geodesics are secondary invariants.  However, their affine parameterizations are not CR invariant.  Instead, if $z$ is an affine parameter for $\gamma$ with respect to $\theta$, then a solution $\hat z=\hat z(z)$ to
\begin{equation}
 \label{eqn:affine_reparameterization}
 \frac{\hat z^{\prime\prime}}{\hat z^\prime} = 2\Upsilon^\prime ,
\end{equation}
is an affine parameter for $\gamma$ with respect to $\htheta$.

Let
\[ \{ p, z \} = \frac{p^{\prime\prime\prime}}{p^\prime} - \frac{3}{2}\left(\frac{p^{\prime\prime}}{p^\prime}\right)^2 \]
denote the Schwarzian derivative of $p(z)$, where primes denote differentiation with respect to $z$.  Recall that $\{p,z\}=0$ if and only if $p=(az+b)/(cz+d)$ for some $ad-bc\not=0$.  Recall also the Chain Rule
\begin{equation}
 \label{eqn:schwarzian_chain_rule}
 \{ p, w \} = \bigl( \{ p, z \} - \{ w, z\} \bigr) \left(\frac{dz}{dw}\right)^{2} .
\end{equation}
We say that $p$ is a \emph{projective parameter} for $\gamma$ if
\begin{equation}
 \label{eqn:projective_parameterization}
 \{ p, z\} = 2i(\gamma^\prime)^\alpha(\gamma^\prime)^\beta A_{\alpha\beta}
\end{equation}
where $z$ is an affine parameter for $\gamma$ with respect to $\theta$ and $A_{\alpha\beta}$ denotes the Tanaka--Webster torsion of $\theta$.  Equations~\eqref{eqn:schwarzian_chain_rule} and~\eqref{eqn:projective_parameterization} imply that the notion of a projective parameter is well-defined independently of the choice of affine parameter and that the projective parameter is uniquely determined up to projective equivalence.

We now show that the projective parameterization $p$ is CR invariant.  To that end, let $\theta$ and $\htheta=e^\Upsilon\theta$ be two contact forms and let $z$ and $\hat z$ be affine parameterizations for $\gamma$ with respect to $\theta$ and $\htheta$, respectively.  By~\eqref{eqn:affine_reparameterization}, it holds that
\[ \frac{\hat z^{\prime\prime}}{\hat z^\prime} = 2\Upsilon^\prime, \]
and hence
\[ \frac{\hat z^{\prime\prime\prime}}{\hat z^\prime} = 2\Upsilon^{\prime\prime} + 4(\Upsilon^\prime)^2 . \]
Since $z$ is an affine coordinate, $\Upsilon^{\prime\prime}=(\gamma^\prime)^\alpha (\gamma^\prime)^\beta\Upsilon_{\alpha\beta}$.  Moreover, we have that
\[ \{ \hat z, z \} = 2\bigl( \Upsilon^{\prime\prime} - (\Upsilon^\prime)^2 \bigr) . \]
Now let $\hat p$ and $p$ be projective parameterizations of $\gamma$ with respect to $\theta$ and $\htheta$, respectively.  Using the Chain Rule~\eqref{eqn:schwarzian_chain_rule}, we deduce that
\begin{align*}
 \left\{ \hat p, p \right\} & = \left( \left\{ \hat p, \hat z \right\} - \left( \left\{ p, z\right\} - \left\{ \hat z, z\right\}\right) \left(\frac{dz}{d\hat z}\right)^2\right) \left(\frac{d\hat z}{dp}\right)^2 \\
 & = 2i\left((\hat\gamma^\prime)^\alpha(\hat\gamma^\prime)^\beta \hA_{\alpha\beta} - \left((\gamma^\prime)^\alpha(\gamma^\prime)^\beta A_{\alpha\beta} + i\left(\Upsilon^{\prime\prime}-(\Upsilon^\prime)^2\right)\right)\left(\frac{dz}{d\hat z}\right)^2\right)\left(\frac{d\hat z}{dp}\right)^2,
\end{align*}
where the second equality uses~\eqref{eqn:affine_reparameterization}.  Recall that the pseudohermitian torsion transforms by~\cite{Lee1988}
\[ \hA_{\alpha\beta} = A_{\alpha\beta} + i\Upsilon_{\alpha\beta} - i\Upsilon_\alpha\Upsilon_\beta . \]
Inserting this into the previous display yields
\[ \{ \hat p, p \} = 0; \]
that is, the projective parameterizations agree.

Finally, note that the above construction is purely local.  In particular, if $\theta$ and $\htheta$ are two pseudo-Einstein contact forms defined in a neighborhood $V$ of a point on a complex null geodesic $\gamma$, then the projective parameterizations of $\gamma\cap V$ determined by $\theta$ and $\htheta$ are projectively equivalent.  These local projective parameterizations combine to give a projective structure on the whole of $\gamma$.

\section{Proof of Theorem~\ref{thm:main}}
\label{sec:details}

We prove Theorem~\ref{thm:main} in the case $n=4$; proofs of the other cases follow analogously.  Let $M=\bR\times\bC^4$ and set
\begin{multline*}
 \theta = \frac{1}{2}dt + \frac{i}{2}\left(z_1\,d\bar z_2 + z_2\,d\bar z_1 + z_3\,d\bar z_4 + z_4\,d\bar z_3 - \bar z_1\,dz_2 - \bar z_2\,dz_1 - \bar z_3\,dz_4 - \bar z_4\,dz_3\right) \\ - i\left(z_1\bar z_1^2\,dz_1 - z_1^2\bar z_1\,d\bar z_1\right) .
\end{multline*}
Let $\{T,Z_1,\dotsc,Z_4\}$ be the basis dual to the admissible coframe $\{\theta,dz^1,\dotsc,dz^4\}$; i.e.\ $T=2\partial_t$ and
\begin{align*}
 Z_1 & := \partial_{z_1} + i\left(\bar z_2 + 2z_1\bar z_1^2\right)\partial_t, \\
 Z_2 & := \partial_{z_2} + i\bar z_1\partial_t, \\
 Z_3 & := \partial_{z_3} + i\bar z_4\partial_t, \\
 Z_4 & := \partial_{z_4} + i\bar z_3\partial_t .
\end{align*}
We define $\bar Z_j$, $j\in\{1,2,3,4\}$, by conjugation; e.g.\ $\bar Z_2 = \delbar_{z_2} - iz_1\partial_t$.  With $T^{1,0}$ the span of $\{Z_1,\dotsc,Z_4\}$ and $T^{0,1}$ its conjugate, we see that $(M,J,\theta)$ is the pseudohermitian manifold $M=r^{-1}(0)$ with contact form $i\partial r$ described in the introduction.

In terms of the coframe $\{\theta,dz^1,\dotsc,dz^4\}$, one readily checks that $\theta$ is a torsion-free contact form and that the only nonvanishing Tanaka--Webster connection one-form $\omega_\alpha{}^\beta$ is
\begin{equation}
 \label{eqn:connection}
 \omega_1{}^2 = 4\bar z_1\,\theta^1 .
\end{equation}
It readily follows that
\begin{equation}
 \label{eqn:curvature}
 R_1{}^2{}_{1\bar 1} = -4
\end{equation}
and all other components of the curvature $R_\alpha{}^\beta{}_{\rho\bar\sigma}$ of $\theta$ vanish.  Therefore the pseudohermitian Ricci tensor $R_{\alpha\bar\beta}:=R_\gamma{}^\gamma{}_{\alpha\bar\beta}$ of $\theta$ vanishes identically.  Hence $\theta$ is pseudo-Einstein and the Chern tensor $S_\alpha{}^\beta{}_{\rho\bar\sigma}$, which is the obstruction to $\theta$ being locally CR flat, agrees with $R_\alpha{}^\beta{}_{\rho\bar\sigma}$.  We thus verify that $(M,\theta)$ is not locally CR flat.

Recall that the Chern tensor $S_\alpha{}^\beta{}_{\rho\bar\sigma}$ is CR invariant.  Consider it as a map $S\colon\Gamma\bigl(T^{1,0}\otimes T^{1,0}\otimes T^{0,1}\bigr)\to\Gamma\bigl(T^{1,0}\bigr)$ by
\[ (U^\alpha, V^\gamma, W^{\bar\sigma}) \mapsto S_\alpha{}^\beta{}_{\gamma\bar\sigma}U^\alpha V^\gamma W^{\bar\sigma} . \]
From~\eqref{eqn:curvature}, we conclude that the image of the Chern tensor of $(M,\theta)$ is the complex span of $Z_2$.  Consider the leaves of the corresponding rank two real distribution $\Real\,\lp Z_2\rp$, which are necessarily complex null curves.  Since $\omega_2{}^\beta=0$ for all $\beta$, we conclude that the leaves are in fact complex null geodesics.  Of these leaves, only $\mF:=\left\{ (0, 0,z_2,0,0) \suchthat z_2\in\bC\setminus\{0\} \right\}$ is invariant under $\Gamma_{\alpha,\beta}$.  Thus $M_{\alpha,\beta}:=M/\Gamma_{\alpha,\beta}$ admits a \emph{unique} closed leaf $\mF_{\alpha,\beta}$ of the real distribution obtained as the real part of the image of the Chern tensor.  Furthermore, since $\theta$ is torsion-free, $z_2$ is the projective parameter of $\mF$.

While $\theta$ is not $\Gamma_{\alpha,\beta}$-invariant, and hence does not descend to $M_{\alpha,\beta}$, it does induce a pseudo-Einstein contact form on any sufficiently small open set in $M_{\alpha,\beta}$.  In particular, $M_{\alpha,\beta}$ locally admits pseudo-Einstein contact forms.  It follows from Section~\ref{sec:markowitz} that $\mF_{\alpha,\beta}$ has a canonical projective structure; indeed, it is the projective structure induced from the projective structure on $\mF$ after taking the quotient by $\Gamma_{\alpha,\beta}$.

Suppose now that $M_{\alpha,\beta}$ and $M_{\tilde\alpha,\tilde\beta}$ are CR equivalent.  Then the uniquely-determined complex null geodesics $\mF_{\alpha,\beta}$ and $\mF_{\tilde\alpha,\tilde\beta}$ are projectively equivalent.  Since the lifts of the respective projective structures to $\mF$ agree, we conclude that $\bC^\ast / \{ z\mapsto e^{3\beta} z\}$ and $\bC^\ast / \{ z\mapsto e^{3\tilde\beta} z\}$ are projectively equivalent.  As $\beta,\tilde\beta<0$, this implies that $\beta=\tilde\beta$.  This completes the proof of Theorem~\ref{thm:main}.

\bibliographystyle{abbrv}
\bibliography{bib}
\end{document}